\documentclass[12pt]{article}
\usepackage{amsfonts,amsmath,amssymb,latexsym}
\newtheorem{theorem}{Theorem}
\newtheorem{proposition}[theorem]{Proposition}
\newtheorem{corollary}[theorem]{Corollary}
\newtheorem{lemma}[theorem]{Lemma}

\newtheorem{remark}[theorem]{Remark}

\newtheorem{examples}[theorem]{Examples}

\newcommand{\R}{\mathbb{R}}

\def\Ie{\mathbb{I}}
\def\dd{\mathrm{d}}

\def\<{\langle}
\def\>{\rangle}
\def\bea{\begin{eqnarray*}}
\def\eea{\end{eqnarray*}}
\def\be{\begin{equation}}
\def\ee{\end{equation}}

\newcommand{\po}{{\hspace*{-1ex}}{\bf .  }}
\def\proof{\noindent{\it Proof: }}
\def\qed{\ifhmode\unskip\nobreak\fi\ifmmode\ifinner\else
\hskip5 pt \fi\fi\hbox{\hskip5 pt \vrule width4 pt height6 pt
depth1.5 pt \hskip 1pt }}

\begin{document}

\title{Conformal Killing graphs with prescribed\\ mean curvature}
\author{M. Dajczer\thanks{Partially supported by
CNPq and Faperj.}\,\,\,\,and\, J. H. de Lira\thanks {Partially supported
by CNPq and Funcap.}}
\date{}
\maketitle

\begin{abstract}
We prove the existence and uniqueness of  graphs with prescribed
mean curvature function in a large class of Riemannian manifolds which comprises spaces endowed with a conformal Killing vector field.
\end{abstract}


\section{Introduction}

Our aim in this paper is to continue the research theme elaborated in
\cite{DHL}  where the existence of a Killing vector field in a manifold permitted to formulate and solve Dirichlet problems associated to the mean curvature equation.  Here, we expand our scope by dealing with conformal Killing vector fields and the corresponding notion of graph. Our main goal is to provide a significant improvement  of the results in \cite{AD} for graphs of constant mean curvature and trivial initial condition. This is achieved by means of a completely different methodology that allows us to generalize, under an unifying perspective, several previous results. 

The pure analytic method  used here enable us to discard or weaken several assumptions   in \cite{AD}. For instance, we allow prescribed mean curvature and nontrivial boundary data. Moreover, we  remove the restriction on the conformal Killing field to be closed.
We also remove the restriction on the Ricci curvature to be minimal in the direction of the Killing field and weaken other requirements on the Ricci tensor. Furthermore, we ruled out some restrictions on the metric of the ambient space so that, in particular, our result applies successfully to a class of product spaces.

To explain the framework of this paper we first fix some
terminology. Let $\bar M^{n+1}$ denote a Riemannian manifold endowed
with a conformal Killing vector field $Y$ whose orthogonal
distribution $\mathcal{D}$ is integrable. Thus, there exists   $\rho\in C^{\infty}(\bar{M})$ such that $\pounds_{Y} \bar
g =2\rho \bar g$, where $\bar g$ is the metric in $\bar M$. It
results that the integral leaves of $\mathcal{D}$ are totally
umbilical hypersurfaces. If in addition $Y$ is closed, then  they are spherical, i.e., have constant mean curvature.

We denote by $\Phi:\Ie\times M^n\to\bar M^{n+1} $ the flow generated
by $Y$, where $\Ie=(-\infty, a)$ is an interval with $a>0$ and $M^n$
is an arbitrarily fixed integral leaf of $\mathcal{D}$
labeled  as $t=0$. It may happen that $a=+\infty$, i.e., the vector
field $Y$ is complete. For instance, this  occurs
when the trajectories of $Y$ are circles and we pass to the
universal cover in an appropriate manner. Since $\Phi_t = \Phi(t,\,\cdot\,)$ is a conformal
map for any fixed $t\in\Ie$, there exists a positive function
$\lambda\in C^\infty(\Ie)$ such that $\lambda(0)=1$ and
$\Phi_t^*\bar g=\lambda^2(t)\bar g$.

Given a bounded domain $\Omega$ in $M$, the \emph{conformal Killing graph}
$\Sigma=\Sigma(z)$ of a function $z$ on $\bar\Omega$ is the hypersurface
$$
\Sigma=\{\bar u=\Phi(z(u),u): u\in\bar\Omega\}.
$$
Proving the existence of a conformal Killing graph with
prescribed mean curvature and boundary requires establishing
apriori estimates. This is accomplished by the use of  Killing cylinders as
barriers. The Killing cylinder $K$ over
$\Gamma=\partial\Omega$ is the hypersurface ruled by the flow lines
of $Y$ through $\Gamma$, that is,
$$
K=\{\bar u=\Phi(t,u):u\in\Gamma\}.
$$

Let $\Omega_0$ denote the largest open subset of points of
$\Omega$ that can be joined to $\Gamma$ by a \emph{unique}
minimizing geodesic. At points of $\Omega_0$,  we denote
$$
\textrm{Ric}_{\bar M}^{rad}(x)=\textrm{Ric}_{\bar M}(\eta,\eta)
$$
where $\textrm{Ric}_{\bar M}$ is the ambient Ricci tensor and $\eta\in T_xM$ is a unit vector tangent to the the unique
minimizing geodesic from $x\in\Omega_0$ to $\Gamma$.

The following result assures the existence of conformal Killing
graphs with prescribed mean curvature $H$ and boundary data $\phi$.
Here, the functions $H$ and $\phi$ are defined on
$\bar\Omega$ and $\Gamma$, respectively. Moreover, $H_{K}$ denotes the mean curvature of  $K$ when calculated pointing inwards.

\begin{theorem}\po\label{main}
Let $\Omega\subset M^n$ be a $C^{2,\alpha}$ bounded domain so that
${\rm Ric}_{\bar M}^{rad}\ge-n\inf_\Gamma H_{K}^2$. Assume  $\lambda_t\geq 0$ and $(\lambda_t/\lambda)_t\geq 0$.
Let $H\in C^\alpha(\Omega)$ and
$\phi\in C^{2,\alpha}(\Gamma)$ be such that $\inf_\Gamma
H_{K}> H\ge 0$ and \mbox{$\phi\leq 0$}.  Then, there exists a unique
function $z\in C^{2,\alpha}(\bar\Omega)$
whose conformal Killing graph has mean curvature function $H$ and
boundary data~$\phi$.
\end{theorem}

Proposition \ref{Hgamma2} below implies that Theorem \ref{main}
holds under weaker but somewhat more technical assumptions on the
ambient Ricci tensor. We also point out that we can prove with minor
modifications   an existence result for functions $H$
depending also on $t$ by imposing the condition $(\lambda H)_t\ge 0$
instead of $\lambda_t\ge 0$. In the case comprised in
Theorem~\ref{main} the condition $\lambda_t\ge 0$ says that the mean
curvatures of the leaves and of the graph have opposite signs.

It is worth to mention that $\bar M$ is conformal to a Riemannian product manifold  $\mathbb{I}\times \tilde{M}$ where $\tilde{M}$ is conformal to $M$. A quite remarkable fact is that the mean curvature equation for a general conformal metric to a product metric like this does not satisfy, in  general, the maximum principle. In fact, in the final remark of the paper we see that the class of metrics which we deal in this paper stands as a borderline for the validity of the elliptic techniques in the treatment of the mean curvature equation.

Assume now that $Y$ is closed or just a Killing field. In each case, a short additional argument to the proof of Theorem \ref{main} yields the same conclusion under weaker assumptions.
For instance,   in the case of a Killing field we prove a generalization of the main result in \cite{DHL} under a weaker assumption on the Ricci curvature.

\begin{corollary}\po\label{main3} Assume that $Y$ is a Killing field. Let $\Omega\subset M^n$ be a $C^{2,\alpha}$ bounded domain such that
${\rm Ric}_{\bar M}^{rad}\ge-n\inf_\Gamma H_{K}^2$.
Let $H\in C^\alpha(\Omega)$ and
$\phi\in C^{2,\alpha}(\Gamma)$ be such that $\inf_\Gamma
H_{K}\geq H\ge 0$.  Then, there exists a unique
function $z\in C^{2,\alpha}(\bar\Omega)$
whose conformal Killing graph has mean curvature function $H$ and
boundary data~$\phi$.
\end{corollary}

The case of closed conformal Killing fields encompasses a broad range of examples, namely,
product and warped ambient spaces, that have been extensively
considered in the recent pertinent literature.
In this case, we can state our result in terms of the Ricci tensor  of $M^n$.

\begin{corollary}\po\label{main2} Assume that $Y$ is a closed conformal Killing field. Let $\Omega\subset M^n$ be a $C^{2,\alpha}$ bounded domain such that $n\,{\rm Ric}_M^{rad}\ge -(n-1)^2\inf_{\Gamma} H_\Gamma^2$. Assume  $\lambda_t\geq 0$ and
$(\lambda_t/\lambda)_t\geq 0$.
Let $H\in C^\alpha(\Omega)$ and
$\phi\in C^{2,\alpha}(\Gamma)$ be such that $\inf_\Gamma
H_{K}> H\ge 0$ and \mbox{$\phi\leq 0$}.  Then, there exists a unique
function $z\in C^{2,\alpha}(\bar\Omega)$
whose conformal Killing graph has mean curvature function $H$ and
boundary data~$\phi$.
\end{corollary}

We finally point out that concerning constant mean curvature Euclidean radial graphs over spherical domains, the above result extends the theorems for minimal radial graphs due to Rado and 
Tausch as well as for constant mean curvature by Serrin and Lopez
(see \cite{Lo} and references therein).


\section{Preliminaries}

Let $(\bar{M}^{n+1},\bar g)$ be a  Riemannian manifold endowed with a
conformal Killing vector field $Y$ whose orthogonal distribution
$\mathcal{D}$ is integrable. Let $\bar\nabla$ denote the Riemannian connection in
$\bar M^{n+1}$  and
$$
\<X,Z\rangle=\bar g(X,Z).
$$
 From $\pounds_{Y}\bar g =2\rho \bar g$ we deduce
the conformal Killing equation \be\label{Lie2}
\<\bar{\nabla}_{X}Y,Z\rangle+\< \bar{\nabla}_{Z}Y,X\>=2\rho
\<X,Z\rangle, \ee where $X,Z\in T\bar M$. It is a standard fact
(cf.\ \cite{Po}) that the conformal factor $\lambda\in C^\infty(\Ie)$ and
$\rho\in C^\infty(\Ie)$ are related by \be\label{hrho}
\rho=\lambda_t/\lambda. \ee Denote
$$
|Y(t,u)|^2=1/\bar\gamma(t,u)\;\;\;\mbox{and}\;\;\;
\gamma(u)=\bar\gamma(0,u).
$$
It follows from (\ref{Lie2}) and (\ref{hrho}) that
\be\label{gamma}
\bar\gamma(t,u)=\gamma(u)/\lambda^2(t).
\ee
We have from (\ref{Lie2}) and the integrability of $\mathcal{D}$ that
\be\label{umbilical}
\<\bar{\nabla}_{X}Y,Z\rangle=\rho\<X,Z\rangle
\ee
for any $X,Z\in\mathcal{D}$. Thus, the leaves $M_t^n=\Phi_t(M)$ are totally umbilical  and the mean curvature $k=k(t,u)$
of $M_t$ with respect to the unit normal vector field $Y/|Y|$ is
\be\label{k}
k=-\frac{\rho}{|Y|}
=-\frac{\lambda_t\sqrt{\gamma}}{\lambda^2}.
\ee

We assign coordinates $x^0=t,x^{1},\ldots,x^{n}$  to points
in $\bar M$ of the form $\bar{u}=\Phi(t,u)$ where
$x^{1},\ldots,x^{n}$ are local coordinates in $M$. Then, the coordinate vector fields are
$$
\partial_0|_{\bar u} =
Y(\bar u)\;\;\;\mbox{and}\;\;\; \partial_i|_{\bar u}=\Phi_{t\,*} \partial_i |_{u}\;\;\mbox{for}\;\;1\le i\le n.
$$
 The components of the ambient line element $\dd s^2$ in terms of these coordinates are
$$
\bar{\sigma}_{00}=\<\partial_0,  \partial_0\> =
|Y|^2=\lambda^2(t)/\gamma(u),\quad \bar{\sigma}_{0i}
=\<\partial_0,\partial_i\rangle= 0,\quad
\bar{\sigma}_{ij}|_{\bar{u}} = \lambda^2(t)\sigma_{ij}|_u,
$$
where $\sigma_{ij}$ are the local components of the metric
$\dd\sigma^2$ in  $M^n$. Setting,
$$
\psi^2(u)= 1/\gamma(u)
$$
we have that $\bar M^{n+1}$ is conformal  to the Riemannian \emph{warped} product manifold $M^n\times_\psi\Ie$
with conformal factor $\lambda$, i.e.,
\be\label{kill-conf}
\dd s^2=\lambda^2(t)(\psi^2(u) \dd t^2 + \dd\sigma^2).
\ee
Finally, after the change of variable
$$
r=r(t)=\int_0^t\lambda(\tau)d\tau,
$$
we have that (\ref{kill-conf}) takes the form of a Riemannian \emph{twisted} product
$$
\dd s^2=\psi^2(u) \dd r^2 + \theta^2(r)\dd\sigma^2
$$
where $\theta(r)=\lambda(t(r))$.
\medskip

The above change of variable is essential to avoid  the coefficients
of the terms of second order in the quasilinear elliptic mean
curvature equation for $\Sigma(z)$ to depend on the function $z$
itself. 

We conclude this sections with a few examples to illustrate
this change of variable.
Observe that all Riemannian manifolds in the examples below fit
the assumptions in Theorem \ref{main}. Moreover, these examples comprise Euclidean and hyperbolic space forms what implies that Theorem \ref{main} assures existence for radial graphs with prescribed mean curvature in these spaces.

\begin{examples}\po {\em Let $\phi\in C^\infty(M)$ be a positive function.\vspace{1ex}

\noindent\text{(a)} By means of the change of variable $e^t=r$,
we have that
$$
\bar{M}= \R_+\times M,\;\;\dd \bar{s}^2=\phi^2(u)\,\dd r^2 + r^2\dd \sigma^2
$$
is isometric to
$$
\tilde{M}= \R\times M,\;\;
\dd \tilde{s}^2=e^{2t}(\phi^2(u)\,\dd t^2 + \dd \sigma^2).
$$
\noindent\text{(b)} By means of the change of variable $t=1-e^{-r}$, we have that
$$
\bar{M}=\R\times M ,\;\;\dd \bar{s}^2
=\phi^2(u)\,\dd r^2 + e^{2r}\dd \sigma^2
$$
is isometric to
$$
\tilde{M}= (-\infty,1)\times M,\;\;
\dd \tilde{s}^2=\frac{1}{(1-t)^2}(\phi^2(u)\,\dd t^2 + \dd \sigma^2).
$$
\noindent\text{(c)} By means of the change of variable $t=c+\log (b\tanh\,(r/2))$,
where $c>0$ and $b^{-1}=\tanh\,(c/2)$, we have that
$$
\bar{M}=\R_+\times M ,\;\;\dd \bar{s}^2
=\phi^2(u)\,\dd r^2 + (\sinh r)^2\dd \sigma^2
$$
is isometric to
$$
\tilde{M}= (-\infty,c+\log b)\times M,\;\;
\dd \tilde{s}^2= (\sinh(2\,\mathrm{argtanh}\,b^{-1}e^{t-c}))^2(\phi^2(u)\,
\dd t^2+ \dd \sigma^2).
$$
}\end{examples}

In the particular case when $Y$ is \emph{closed} conformal Killing field, that is, when
$$
\<\bar{\nabla}_{X}Y,Z\rangle=\rho \<X,Z\rangle,
$$
we have that $\gamma$ is constant.  Thus, in this case, $\bar M^{n+1}$ has a warped
product structure~and is~conformal with conformal factor $\lambda$
to a Riemannian product manifold  $\Ie\times M^n$.
Observe that now  the leaves of $\mathcal{D}$ are spherical, that is, totally umbilical with constant mean curvature $k=k(t)$.

\section{Killing cylinders}

Let $\Omega\subset M^n$ be a bounded domain  with regular
boundary~$\Gamma$. The \emph{Killing cylinder} $K$ over $\Gamma$
determined by the conformal Killing field $Y$ is the hypersurface
defined by
$$
K=\{\Phi(t,u): t\in\Ie,\; u\in \Gamma\}.
$$

Let $t^1,\ldots, t^{n-1}$ be local coordinates for $\Gamma$. We
denote by $(\tau_{ij})$  the components of the metric in $\Gamma$
with respect to these coordinates. It results that $t,t^1,\ldots,
t^{n-1}$ are local coordinates for $K$. Let $\eta$ be the
inward unit normal vector along $\Gamma$ as a submanifold of $M$. Then
$$
\bar\eta=\frac{1}{\lambda}\Phi_{t_*}\eta
$$
is an unit normal vector field to $K$. Thus,
$$
\<\bar\eta,\partial_t\>=0=\<\bar\eta,\partial_{t^i}\>
$$
where $\partial_t=\partial/\partial t$ and
$\partial_{t^i}=\partial/\partial t^i$. We deduce from
(\ref{umbilical}) that
$$
\<\bar\nabla_{\partial_{t^i}}\partial_t,\bar\eta\>
=\rho\<\partial_{t^i},\bar\eta \> =0.
$$
Hence $\partial_t$ is a principal  direction of $K$ with
corresponding principal curvature
\be\label{kappa}
\kappa=\bar\gamma\<\bar\nabla_YY,\bar\eta\>.
\ee
It follows from (\ref{Lie2}) and (\ref{gamma}) that
\be\label{kappa2}
\kappa=-\frac{1}{2}\bar\gamma\bar\eta(\bar\gamma^{-1})
=-\frac{1}{2}\gamma\bar\eta(\gamma^{-1})
=\frac{1}{2\gamma}\bar\eta(\gamma)=
\frac{1}{\lambda}\eta(\log\sqrt{\gamma}).
\ee

It was shown in \cite{LN} that  the distance function
$d(u)=\textrm{dist}(u,\Gamma)$ in $\Omega_0$ has the same regularity
as $\Gamma$.  Hence, local coordinates in $\bar M$ near $K$ can be
defined setting $t^0=t$ and $t^n=d$. We denote by $(t_{ij})$ the
components of the metric for these coordinates. Thus,
$$t_{ij}(t,u)=\lambda^2(t)\tau_{ij}(u)\;\;\mbox{for}\;\; 1\leq i,j\leq n-1.
$$

\begin{lemma}\label{Hgamma}\po
The mean curvature $H_K$ of the Killing cylinder $K$ is given by
\be\label{mean}
nH_K(t,u)=\kappa(t,u) +
\frac{n-1}{\lambda(t)}H_{\Gamma}(u).
\ee
\end{lemma}

\proof  We have that
\begin{eqnarray}
 nH_K\!\!\!&= &\!\!\! \kappa
+ t^{ij}\<\bar\nabla_{\partial_{t^i}}\partial_{t^j},
\bar\eta\>|_{(t,u)}=\kappa + \lambda^{-2}\tau^{ij}
\<\Phi_{t*}\nabla_{\partial_{t^i}}\partial_{t^j}|_u,
\lambda^{-1}\Phi_{t*}\eta|_u\>\nonumber\\
\!\!\!&= &\!\!\!  \kappa +
\lambda^{-1}\tau^{ij}\<\nabla_{\partial_{t^i}}
\partial_{t^j}|_u,\eta|_u\>,\nonumber
\end{eqnarray}
and the proof follows.\qed\vspace{1,5ex}

We denote by $\Gamma_\epsilon$ and $K_\epsilon$ the level sets
$d=\epsilon$ in $M$ and $\bar M$, respectively.
By $H_{K_\epsilon}$ we denote the mean curvature of the Killing cylinder $K_\epsilon$ over $\Gamma_\epsilon$.

\begin{proposition}\label{Hgamma2}\po
Assume that the Ricci curvature tensor of $\bar M$ satisfies
\be\label{ric-amb} \inf_{\Omega_0}\big\{\emph{Ric}_{\bar M}^{rad}
+\big(nk^2-\sqrt{\gamma}\,k_t\big)|_{t=0}\big\} \ge -n\inf_{\Gamma}
H_K^2. \ee Then, $H_{K_\epsilon}|_{x_0}\geq H_K|_{y_0}$ if $y_0\in
\Gamma$ is the closest point to a given point
$x_0\in\Gamma_\epsilon\subset \Omega_0$.
\end{proposition}

\proof The coordinate $d$-curve in
$\Phi(t,u)$ is the image by $\Phi_t$ of the coordinate $d$-curve
passing through $u\in M$. Thus,
$$
\bar\eta|_{\Phi(t,u)}=\frac{1}{\lambda}\Phi_{t\,*}(u)\partial_d|_u
=\frac{1}{\lambda}\partial_d|_{\Phi(t,u)}
=\frac{1}{\lambda}\partial_{t^n}|_{\Phi(t,u)}.
$$

Extend $\bar\eta$ near $K$ as the velocity vector field of the geodesics in $M_t$
departing orthogonally from $K\cap M_t$. We obtain for $1\le i,j\le n-1$ that
\begin{eqnarray*}
\lambda^2\bar\eta\<\bar\nabla_{\partial_{t^i}}\bar\eta,
\partial_{t^j}\>  \!\!\!&=&\!\!\!
\partial_{t^n} \<\bar\nabla_{\partial_{t^i}}
\partial_{t^n}, \partial_{t^j}\>\\
\!\!&= &\!\!\<\bar\nabla_{\partial_{t^i}}
\bar\nabla_{\partial_{t^n}}
\partial_{t^n}, \partial_{t^j}\>
+\<\bar R(\partial_{t^n},
\partial_{t^i})\partial_{t^n}, \partial_{t^j}\>+ \<
\bar\nabla_{\partial_{t^i}}\partial_{t^n}, \bar\nabla_{\partial_{t^j}}\partial_{t^n}\>\\
\!\!\!&=&\!\!\! -\lambda^2(\bar\gamma \<\bar\nabla_{\bar\eta}\bar\eta,Y\>
\<\bar\nabla_{\partial_{t^i}}\partial_{t^j}, Y\>
+\<\bar R(\partial_{t^i},\bar\eta)\bar\eta,\partial_{t^j}\>
-\<\bar\nabla_{\partial_{t^i}}\bar\eta,
\bar\nabla_{\partial_{t^j}}\bar\eta\>).
\end{eqnarray*}
Using  (\ref{umbilical}) and (\ref{k}), we have
\be\label{first}
\bar\eta \<\bar\nabla_{\partial_{t^i}}\bar\eta,\partial_{t^j}\>
=-k^2t_{ij}
-\<\bar R( \partial_{t^i},\bar\eta)\bar\eta,\partial_{t^j}\>
+ \<A^2_\epsilon\partial_{t^i}, \partial_{t^j}\>
\ee
where $A_\epsilon$ denotes the Weingarten map of $K_\epsilon$ relative to $\bar\eta$.

For the remaining  case  $i=j=0$, we have
$$
\bar\eta \<\bar\nabla_{\partial_{t^0}}\bar\eta,\partial_{t^0}\> =
\<\bar\nabla_{Y}\bar\nabla_{\bar\eta}\bar\eta, Y\>+\<\bar
R(\bar\eta,Y)\bar\eta,Y\>
+\<\bar\nabla_{[\bar\eta,Y]}\bar\eta,Y\rangle
+\<\bar\nabla_Y\bar\eta,\bar\nabla_{\bar\eta}Y\rangle.
$$
However,
\be\label{corch}
[\bar\eta,Y]=-[Y,\bar\eta]=-[\partial_{t^0},
\lambda^{-1}\partial_{t^n}]
=\frac{\lambda_t}{\lambda^2}\partial_{t^n}=\rho\bar\eta.
\ee
Thus,
$$
\<\bar\nabla_{[\bar\eta,Y]}\bar\eta,Y\rangle
=\rho\langle\bar\nabla_{\bar\eta}\bar\eta, Y\rangle=-\rho^2.
$$
Using (\ref{corch}) we have
$$
\<\bar\nabla_Y\bar\eta,\bar\nabla_{\bar\eta}Y\rangle
=\<\bar\nabla_Y\bar\eta,\bar\nabla_Y
{\bar\eta}\rangle+\<\bar\nabla_Y\bar\eta,[{\bar\eta},Y]\rangle
=\langle A_\epsilon^2Y,Y\rangle
$$
and
$$
\<\bar\nabla_{Y}\bar\nabla_{\bar\eta}\bar\eta, Y\>
= Y\<\bar\nabla_{\bar\eta}\bar\eta, Y\>
-\bar\gamma\<\bar\nabla_{\bar\eta}\bar\eta,Y\>\<\bar\nabla_{Y}Y,Y\>
=-Y(\rho)+\rho^2.
$$
It follows that
\be\label{second}
\bar\eta \<\bar\nabla_{\partial_{t^0}}\bar\eta,
\partial_{t^0}\>  =-Y(\rho)-\<\bar R(Y,\bar\eta)\bar\eta,Y\>+\<A^2_\epsilon Y,Y\>.
\ee
We also have
$$
\bar\eta \<\bar\nabla_{\partial_{t^i}}\bar\eta,\partial_{t^j}\>
=-\<\bar\nabla_{\bar\eta} A_\epsilon \partial_{t^i},\partial_{t^j}\>
- \<A_\epsilon\partial_{t^i}, \bar\nabla_{\bar\eta}\partial_{t^j}\>
= -\<(\bar\nabla_{\bar\eta} A_\epsilon)\partial_{t^i},
\partial_{t^j}\> +2\<A^2_\epsilon\partial_{t^i},\partial_{t^j}\>
$$
for $0\le i,j\le n-1$. Taking  traces with respect to the metric $(t_{ij})$ in $K$
with $t^{00}=\bar\gamma$ and $t^{0i} =0$,  and using
(\ref{first}) and (\ref{second}) gives
\begin{eqnarray*}
\textrm{tr}\bar\nabla_{\bar\eta} A_\epsilon \!\!\!&=&\!\!\!  t^{ij}
\<(\bar\nabla_{\bar\eta} A_\epsilon)\partial_{t^i},\partial_{t^j}\>
=-t^{ij}\bar\eta \<\bar\nabla_{\partial_{t^i}}
\bar\eta,\partial_{t^j}\>
+2t^{ij}\<A^2_\epsilon\partial_{t^i},\partial_{t^j}\>\\
\!\!\!&=&\!\!\!
-\bar\gamma\bar\eta\<\bar\nabla_{\partial_{t^0}}\bar\eta,
\partial_{t^0}\>
+2\bar\gamma\<A^2_\epsilon\partial_{t^0},\partial_{t^0}\>
- \sum_{i,j=1}^{n-1}t^{ij}(\bar\eta \<
\bar\nabla_{\partial_{t^i}}\bar\eta,\partial_{t^j}\>
-2\<A^2_\epsilon\partial_{t^i},\partial_{t^j}\>)\\
\!\!\!&=&\!\!\! \bar\gamma Y(\rho) +(n-1)k^2+\textrm{Ric}_{\bar
M}(\eta,\eta)+\textrm{tr}A^2_\epsilon.
\end{eqnarray*}
However,
\begin{equation}
\label{kt} \sqrt{\bar\gamma}\,k_t=\frac{Y(k)}{|Y|}
=-\frac{1}{|Y|^2}Y(\rho)+k
Y\Big(\frac{1}{|Y|}\Big)
 = -\bar\gamma Y(\rho)+k^2.
\end{equation}
We conclude that
$$
\textrm{tr}\bar\nabla_{\bar\eta} A_\epsilon
=\textrm{tr}A^2_\epsilon
+{\rm Ric}_{\bar M}(\bar\eta,\bar\eta)
+n k^2-\sqrt{\bar\gamma}\,k_t.
$$
Since $n\dot H_{K_\epsilon}
=\bar\nabla_{\bar\eta}\textrm{tr}A_\epsilon
=\textrm{tr}\bar\nabla_{\bar\eta} A_\epsilon$, at $d=\epsilon$, it follows that
$$
n\dot H_{K_\epsilon}\ge nH^2_{K_\epsilon}
+{\rm Ric}_{\bar M}(\bar\eta,\bar\eta)+ nk^2-\sqrt{\bar\gamma}\,k_t.
$$
Therefore, at $t=0$ and by the assumption on $\textrm{Ric}_{\bar M}^{rad}$
there exist constants $c>0$ and $d_0>0$ such that
$$
\dot H_{K_d}\ge c(H_{K_d}-\inf_{\Gamma} H_{K})
$$
for $d\in [0,d_0]$. Hence, $H_{K_d}$ does not decrease with
increasing $d$. \qed

\section{Conformal Killing graphs}

The \emph{conformal Killing graph} $\Sigma^n$ of a function
$z\colon\,\bar\Omega\subset M^n \to\Ie$  is the hypersurface in
$\bar M^{n+1}$ given by
$$
\Sigma^n=\{X(u)=\Phi(z(u),u):\ u\in \bar\Omega \}.
$$
We show next that the partial differential equation for
a prescribed mean curvature function $H$ in  $\bar\Omega$ is the
quasilinear elliptic equation of divergence form
\be\label{e11}
{\mathrm{div}}\, \bigg(\frac{\nabla z}{\sqrt{\gamma+|\nabla
z|^2}}\bigg) -\frac{1}{\sqrt{\gamma+|\nabla z|^2}}
\bigg(\frac{\<\nabla\gamma,\nabla z\>}
{2\gamma}+n\gamma\rho\bigg)-n\lambda H=0.
\ee
Recall that $\rho=\lambda_t/\lambda$.
Here, the gradient $\nabla$ and the divergence $\textrm{div}$ are
differential operators in the
leaf $M^n$ endowed with the metric $d\sigma^2$.

A sufficient condition to have a maximum principle for (\ref{e11})
(see Theorem 10.1 in \cite{GT}) is
\be\label{cond-1}
(\lambda_t/\lambda)_t=\rho_t\ge
0\;\;\;\mbox{and}\;\;\; \lambda_t H \ge 0.
\ee

The graph $\Sigma$ is parametrized in terms of local coordinates by
$$
X(u)\in \Sigma \mapsto \big(z(x^1,\ldots,x^n),x^1,\ldots,x^n\big).
$$
Thus, the tangent space to $\Sigma$ at $X(u)$
is spanned by the vectors
\be\label{vector}
X_i=z_i\partial_0|_{X(u)}+\partial_i|_{X(u)},
\ee
where $z_i=\partial z/\partial x^i$. Hence,
the metric induced on $\Sigma$ is
$$
g_{ij}|_{X(u)}=
\lambda^2(z(u))\bigg(\sigma_{ij}(u)
+\frac{z_iz_j}{\gamma(u)}\bigg).
$$
The inverse is
$$
g^{ij}|_{X(u)} = \frac{1}{\lambda^2(z(u))}
\bigg(\sigma^{ij}(u)-\frac{z^iz^j}{\gamma(u)
+|\nabla z|^2}\bigg)
$$
where $z^i=\sigma^{ik}z_k$ and $|\nabla z|^2=z^jz_j$ as usual.

Fix the orientation on $\Sigma$ given by the unit normal vector
field \be\label{N} N|_{X(u)}=\frac{1}{\lambda^2W}\,
\big(\gamma(u)\partial_0|_{X(u)}-\Phi_{z(u)\,*} \nabla z(u)\big),
\ee where
$$
\lambda^2W^2 = \gamma + |\nabla z|^2.
$$
Notice that  $\<N,Y\>>0$.
We compute the second fundamental form
$$
a_{ij}=\<\bar\nabla_{X_i}X_j, N\>
$$
of $\Sigma$. From (\ref{vector}) and
(\ref{N}) we obtain
\begin{eqnarray*}
\lambda^2Wa_{ij}= \hspace*{-5ex}&{} &\gamma z_{ij}\<\partial_0,\partial_0\>
+\gamma z_iz_j\<\bar\nabla_{\partial_0}\partial_0,\partial_0\>
+\gamma z_j\<\bar\nabla_{\partial_i}\partial_0,\partial_0\>
+\gamma z_i\<\bar\nabla_{\partial_j}\partial_0,\partial_0\> \\
\!\!\! &+&\!\!\!\gamma
\<\bar\nabla_{\partial_i}\partial_j,\partial_0\>
-z_iz_j\<\bar\nabla_{\partial_0}\partial_0,\Phi_{z(u)*} \nabla
z\>-z_j\<\bar\nabla_{\partial_i}\partial_0,
\Phi_{z(u)*}\nabla z\>\\
\!\!\!&-&\!\!\! z_i\<\bar\nabla_{\partial_j}\partial_0,
\Phi_{z(u)*}\nabla z\>
-\<\bar\nabla_{\partial_i}\partial_j,\Phi_{z(u)*}\nabla z\>.
\end{eqnarray*}
The Levi-Civita connections in  $M$ and $M_t$ are determined by
the same Christoffel symbols since
$\bar \sigma_{ij}|_{\bar u}=\lambda^2\sigma_{ij}|_u$
if  $u\in M$ and $\bar u =\Phi_t(u)\in M_t$.
We have from (\ref{hrho}) and (\ref{umbilical}) that
$$
\<\bar\nabla_{\partial_i}\partial_j,\partial_0\>
=-\rho(z(u))\<\partial_i|_{X(u)},\partial_j|_{X(u)}\>
= -\lambda\lambda_t(z(u))\sigma_{ij}(u).
$$
The terms involving the flow lines acceleration are
$$
\<\bar\nabla_{\partial_0}\partial_0,\partial_0\>|_{X(u)}
=\frac{1}{2}\partial_t|_{t=z(u)}(\lambda^2/\gamma)
= \frac{\lambda\lambda_t(z(u))}{\gamma(u)}
$$
and
$$
\<\bar\nabla_{\partial_i}\partial_0,\partial_0\>|_{X(u)}
=\frac{1}{2}\partial_i|_{X(u)}
\big(\lambda^2/\gamma\big)
 = -\frac{\lambda^2(z(u))\gamma_i(u)}{2\gamma^2(u)}.
$$
Similarly,
$$
\<\bar\nabla_{\partial_i}\partial_j|_{X(u)},\Phi_{z(u)\,*}
\nabla z(u)\> =\<\Phi_{z(u)\,*}
\bar\nabla_{\partial_i}\partial_j|_{u},\Phi_{z(u)\,*}
\nabla z(u)\>=
\lambda^2(z(u))\<\nabla_{\partial_i}\partial_j|_u, \nabla z|_u\>
$$
and
$$
\<\bar\nabla_{\partial_i}\partial_0,\Phi_{z(u)*}\nabla z\>
=\<
\bar\nabla_{\partial_i}\partial_0|_{X(u)}, z^j \partial_j|_{X(u)}\>
= z^j \rho\<\partial_i|_{X(u)},\partial_j|_{X(u)}\>
=z_i\lambda\lambda_t(z(u)).
$$
Since $z_{i;j}=z_{ij}-\<\nabla_{\partial_i}\partial_j,\nabla z\>$
are the Hessian components, we have
$$
Wa_{ij}= z_{i;j}-\frac{\lambda_t}{\lambda}z_iz_j
-\frac{\lambda_t}{\lambda}\gamma\sigma_{ij}
-\frac{\gamma_i}{2\gamma}z_j-\frac{\gamma_j}{2\gamma}z_i
-\frac{\gamma_k}{2\gamma^2}z_iz_jz^k.
$$
We easily obtain
\begin{eqnarray}\label{old}
\lambda^4W^{3}a^{i}_{k}
\!\!\!&=&\!\!\!\lambda^4W^{3}g^{ij}a_{jk}\nonumber\\
\!\!\!&=&\!\!\!\big(\lambda^2W^2\sigma^{ij}-z^iz^j\big)z_{j;k}
-\frac{1}{2}z^i\gamma_k
-\lambda^2W^2\left(\frac{\gamma^i}{2\gamma}z_k
+\gamma\rho\delta^i_k\right).
\end{eqnarray}
Taking traces after dividing both sides by $\lambda^3W^3$ yields
\be\label{basic2}
n\lambda H=\frac{1}{\lambda
W}\left(\sigma^{ij}-\frac{z^iz^j}{\lambda^2W^2}\right)z_{j;i}
-\frac{\gamma_iz^i}{2\lambda^3W^3} -\frac{1}{\lambda
W}\left(\frac{\gamma_iz^i}{2\gamma}
+n\gamma\rho\right).
\ee
Equivalently, we have
\be\label{useful}
\mathcal{Q}[z]:=
\bigg(\frac{z^i}{\sqrt{\gamma+z^kz_k}}\bigg)_{;i}
-\frac{1}{\sqrt{\gamma+z^kz_k}}
\bigg(\frac{\gamma_iz^i}{2\gamma}+n\gamma\rho\bigg)-n\lambda H=0
\ee
as we wished.

\section{The proof of Theorem 1}

Finding a conformal Killing graph $\Sigma(z)$ with prescribed mean curvature function $H$ and boundary data $\phi$ amounts to solve
the Dirichlet problem
\be\label{dirichlet}
\left\{\begin{array}{l}
\!\mathcal{Q}[z] = 0
\vspace{1ex}\\
\!z|_\Gamma=\phi.
\end{array} \right.
\ee

\noindent \emph{Proof of Theorem \ref{main}:}
We apply the well-known continuity method to the family parametrized by $\tau\in[0,1]$ of Dirichlet problems
\be\label{dirichlet2}
\left\{\begin{array}{l}
\!\mathcal{Q}_\tau [z]=0,
\vspace{1ex}\\
\!z|_\Gamma=\tau\phi
\end{array} \right.
\ee
where
$$
\mathcal{Q}_\tau [z] = {\mathrm{div}}\, \bigg(\frac{\nabla z}{\sqrt{\gamma+|\nabla
z|^2}}\bigg) -
\frac{\<\nabla\gamma,\nabla z\>}
{2\gamma\sqrt{\gamma+|\nabla z|^2}}-\tau\bigg(\frac{n\gamma\rho}{\sqrt{\gamma+|\nabla z|^2}}+n\lambda H\bigg).
$$

 Let $I$ be the subset of $[0,1]$ consisting of values
of $\tau$ for which the Dirichlet problem (\ref{dirichlet2}) has a $C^{2,\alpha}$
solution.  Then, the proof  reduces to show that $I=[0,1]$. First,
observe that $I$
is nonempty since $z=0$ is a solution for $\tau=0$. Moreover, we have that $I$ is open in view of (\ref{cond-1}). The difficult part is
to show that $I$ is closed. This follows from the a priori estimates given below in Propositions \ref{C0}, \ref{C1b} and \ref{C1i}
and standard theory of quasilinear elliptic partial differential equations \cite{GT}.\qed
\begin{remark}\po {\em We point out that our existence results still hold if $\phi$ is only
assumed continuous. We may approximate $\phi$ uniformly by smooth
boundary data and use the interior gradient estimate to obtain
strong convergence on compact subsets of~$\Omega$. A local barrier
argument shows that the limiting solutions achieves the given
boundary data.
}\end{remark}

\section{Height estimates}

In this section, we obtain an a priori height estimate.

\begin{proposition}\po\label{C0} Under the assumptions of
Theorem \ref{main} there exists a positive constant $C=C(\Omega, H)$ independent of $\tau$ such that
$$
|z_\tau|_0 \le C+|\phi|_0
$$
for any solution $z_\tau$ of the Dirichlet problem (\ref{dirichlet2}).
\end{proposition}

\noindent \emph{Proof:} In view of (\ref{cond-1}) we may apply the
Comparison Principle (cf.\ Theorem  10.1 in \cite{GT}) to  (\ref{dirichlet2}). From
$$
\mathcal{Q}_\tau[\tau\sup\phi]\le 0 = \mathcal{Q}_\tau[z_\tau]
\;\;\;\mbox{and}\;\;\; z_\tau|_{\Gamma}\le \tau\sup\phi,
$$
we conclude that $z_\tau\le \tau\sup\phi$.
Similarly, for a solution $z$ of (\ref{dirichlet}) we obtain from
$$
\mathcal{Q}_\tau[z]\ge 0=\mathcal{Q}_\tau[z_\tau]\;\;\;\mbox{and}\;\;\;
z|_{\Gamma}\le z_\tau|_{\Gamma}
$$
that $z\le z_\tau$.
\medskip

Next we construct barriers on $\Omega_0$ which are subsolutions to
(\ref{dirichlet}) of the form
\be\label{function}
\varphi(u)=\inf_\Gamma \phi + f(d(u))
\ee
where $d(u)=\textrm{dist}(u,\Gamma)$ and the real function $f$ will be chosen later.
Hence,
\be\label{deri}
\varphi_i = f' d_i \;\;\;\mbox{and}\;\;\;
\varphi_{i;j}=f'' d_i d_j + f' d_{i;j}.
\ee
At points in $\Omega_0$, we have
\be\label{two}
|\nabla d|=1.
\ee
It follows that
\be\label{three}
d^id_{i;j}=0
\ee
and
$$
2\<\nabla_{\partial_d}\nabla d,
\partial_d\>=\partial_d |\nabla d|^2=0.
$$
Moreover,
\be\label{H}
d^{i}_{\textrm{ };i}=\sigma^{ij}d_{i;j}
=\sigma^{ij}\<\nabla_{\partial_i}\nabla d,\partial_j\>
=-(n-1)H_{\Gamma_\epsilon},
\ee
where $H_{\Gamma_\epsilon}$ denotes the mean curvature of
$\Gamma_\epsilon\subset \Omega_0$ with respect to $\eta$.

Combining
$$
\<\nabla\gamma, \nabla z\>=
-\frac{2}{|Y|^4}\<\bar\nabla_{\nabla z} Y, Y\>
=2\gamma^2\<\bar\nabla_Y Y,\nabla z\rangle.
$$
with (\ref{basic2}) yields
\be\label{eq}
\mathcal{Q}[\varphi]=
\frac{1}{U}\Big(\varphi^{i}_{\textrm{ };i} -\frac{\varphi^i
\varphi^j\varphi_{i;j}}{U^2}\Big)
-\frac{\gamma}{U^3}(\gamma+U^2)\<
\bar\nabla_{Y}Y, \nabla
\varphi\>-\frac{n\gamma\rho}{U}-n\lambda H.
\ee
where
$$
U=\lambda W=\sqrt{\gamma+ f'^2}.
$$
Using (\ref{deri}) and (\ref{H}) we obtain
$$
\mathcal{Q}[\varphi]
=\frac{\gamma}{U^3}(f''-\gamma\<\bar\nabla_{Y}Y,\eta\>f')-
\frac{f'}{U}((n-1)H_{\Gamma_\epsilon} +\gamma\<
\bar\nabla_{Y}Y,\eta\>)-\frac{n\gamma\rho}{U}-n\lambda H.
$$

We chose in (\ref{function}) the test function
$$
f=\frac{e^{DB}}{D}\big(e^{-Dd}-1\big)
$$
where $B>\textrm{diam}(\Omega)$ and $D>0$ is a constant to be
chosen later. Then,
$$
f'=-e^{D(B-d)}\;\;\;\mbox{and}\;\;\; f''=-Df'.
$$
Using (\ref{kt}) and that $\rho_t\geq 0$ by assumption, we obtain
$$
nk^2-\sqrt{\bar\gamma}k_t= (n-1)k^2+\bar\gamma \rho_t\geq 0.
$$
It follows from Proposition \ref{Hgamma2} that
$$
H_{K_\epsilon}\ge H^*> H\geq 0
$$
where $H^*=\inf_\Gamma H_{K}>0$. Since  $\lambda(\varphi)\leq 1$, we obtain using Lemma \ref{Hgamma} that
$$
\mathcal{Q}[\varphi]\geq -\frac{\gamma f'}{U^3}(D+\kappa_\epsilon)
-\frac{f'}{U}nH^*-\frac{n\gamma\rho}{U}-nH.
$$
We require $D>\sup_{\Omega_0} |\kappa_\epsilon|$ and denote
$nD_0=D+\kappa_\epsilon$. Thus $\mathcal{Q}[\varphi]>0$ if
\be\label{dif} HU^3 < -H^*f'U^2-\gamma\rho U^2-\gamma D_0 f'.
\ee
Since $f'\to -\infty$ as $D\to +\infty$,  we conclude that for $D$
sufficiently large the inequality holds.
Hence, we have shown that
$$
\mathcal{Q}[\varphi]>0=\mathcal{Q}[z]\;\;\;\mbox{and}\;\;\;\varphi|_{\Gamma}\leq
z|_{\Gamma}
$$

To prove that  $\varphi\le z$ on $\bar\Omega$ we just
follow the reasoning in the proof of Lemma~6 in \cite{DHL}  (see
\cite{GT}, p.\ 171). We conclude that $\varphi$ is a continuous
subsolution for the Dirichlet problem (\ref{dirichlet}).\qed

\section{Boundary gradient estimates }

In this section, we establish an a priori gradient estimate along
the boundary of the domain.

\begin{proposition}\po\label{C1b}
Under the assumptions of
Theorem \ref{main} there exists a positive constant $C=C(\Omega,
H,\phi, |z|_0)$ independent of $\tau$ such that
$$
\sup_\Gamma |\nabla z_\tau|\le C
$$
for any solution $z_\tau$ of the Dirichlet problem (\ref{dirichlet2}).
\end{proposition}

\proof We argue for $\tau=1$. We use
barriers of the form $w+\phi$ along a tubular neighborhood
$\Omega_\epsilon$ of $\Gamma$ where  we extended  $\phi$ to $\Omega_\epsilon$ by taking $\phi(t^i,d)=\phi(t^i)$.
We choose $w=f(d)$ where
$$
f(d) = -\tilde{\mu}\ln (1+\mu d)
$$
and $\mu>0$, $\tilde{\mu}>0$ are constants. Hence,
$$
f'  = \frac{-\mu\tilde{\mu}}{1+\mu d}\quad\textrm{ and
}\quad f''=\frac{1}{\tilde{\mu}}f'^2.
$$
We choose $\tilde{\mu} = c/\ln (1+\mu)$ with $c>0$ to be chosen later. Hence,
$$
f(\epsilon)=\frac{-c\ln (1+\epsilon\mu)}{\ln(1+\mu)}\to -c \;\;\;\mbox{as}\;\;\;\mu\to +\infty\;\;\mbox{for}\;\;\epsilon>0,
$$
and
\be\label{infty}
 f'(0)=\frac{-c\mu}{\ln(1+\mu)} \to -\infty\;\;\;\mbox{as}\;\;\;\mu\to +\infty.
\ee

A simple estimate using (\ref{eq}) gives
\begin{eqnarray*}
\mathcal{Q}[w+\phi]\!\!\!&=&\!\!\! a^{ij}(x,\nabla w+\nabla
\phi)(w_{i;j}+\phi_{i;j})
+b(x,\nabla w+\nabla\phi)-n\lambda H\nonumber\\
\!\!&\ge &\!\! a^{ij}w_{i;j}+\Lambda|\phi|_{2,\alpha} +b-n\lambda H.
\end{eqnarray*}
Here $\Lambda=\gamma/U^3$ is the lowest eigenvalue of the matrix
$$
a^{ij} =\frac{\delta^{ij}}{U}
-\frac{1}{U^3}(w^i+\phi^i)(w^j+\phi^j)
$$
and
$$
U^2=\theta+f'^2\;\;\mbox{where}\;\;\theta=\gamma+|\nabla\phi|^2
$$
from (\ref{two}). Moreover,
$$
b =-\frac{\gamma}{U^3}\big(\gamma+U^2\big)\<\bar\nabla_Y Y, \nabla
w + \nabla\phi\> -\frac{n\gamma\rho}{U}.
$$
It follows from (\ref{two}) and (\ref{three})  that
$$
w^iw^jw_{i;j}=f'^2 d^id^j(f'' d_i  d_j + f'd_{i;j})
=f'^2f'',
$$
$$
w^i\phi^jw_{i;j}=f'f''\<\nabla d,\nabla \phi\>=0
$$
and
$$
\phi^i\phi^j w_{i;j}=f'\phi^i\phi^jd_{i;j}.
$$
Since
$$
\Delta w = f'' +f'\Delta d = f'' -(n-1)f'
H_{\Gamma_d},
$$
we obtain
$$
a^{ij}w_{i;j}= -(n-1)\frac{f'}{U}H_{\Gamma_d} +
\frac{f''}{U^3}(\gamma+|\nabla\phi|^2) -\frac{f'}{U^3}
\phi^i\phi^jd_{i;j}.
$$
Moreover, from (\ref{kappa}) and $\nabla w=f'\eta$ we have
$$
b= -\frac{f'}{U}\Big(\frac{\gamma}{U^2}+1\Big)\kappa
-\frac{\gamma}{U}\Big(\frac{\gamma}{U^2}
+1\Big)\<\bar\nabla_Y Y,\nabla \phi\> -\frac{n\gamma\rho}{U}.
$$
Using Lemma \ref{Hgamma}, we obtain
$$
\begin{array}{l}
\mathcal{Q}[w+\phi]U^3\ge  \displaystyle{ -n(f'H_{K_d} +\lambda
HU)U^2
+\gamma|\phi|_{2,\alpha}-n\gamma\rho U^2}-\gamma\<\bar\nabla_{Y}Y,\nabla \phi\> U^2\vspace*{1.5ex}\\
\hspace*{7ex} \displaystyle{
+f''(\gamma+|\nabla\phi|^2 )-f'\gamma \kappa -f'\phi^i\phi^j d_{i;j}
-\gamma^2\<\bar\nabla_{Y}Y,\nabla\phi\>}.
\end{array}
$$
Since $\phi\le 0$, we have $\lambda(\phi)\le 1$ and
$\rho(\phi)\leq\rho_0$.   At points of $\Gamma$, we obtain
$$
\begin{array}{l}
\mathcal{Q}[w+\phi]U^3 \ge
\displaystyle{-n(f'H_K+H\sqrt{\theta+f'^2})(\theta+f'^2)
+\frac{1}{c}\ln (1+\mu)(\gamma+|\nabla\phi|^2)f'^2}
\vspace*{1.5ex}\\
-\gamma(n\rho_0
+\<\bar\nabla_YY,\nabla\phi\>)(\theta+f'^2) -(\gamma
\kappa +\phi^i\phi^j d_{i;j})f'+\gamma|\phi|_{2,\alpha}
-\gamma^2\<\bar\nabla_{Y}Y,\nabla\phi\>
\end{array}
$$
where $f'=f'(0)$ satisfies (\ref{infty}). It is easy to see using $\inf_\Gamma H_{K}\geq H\geq 0$ that choosing $\mu$ large enough and then $c$ large enough assures that $\mathcal{Q}[w+\phi]>0$ on a small tubular neighborhood $\Omega_\epsilon$ of $\Gamma$ and that $w+\phi\le z$ on both boundary components. Therefore, $w+\phi$ is a locally defined lower barrier for the Dirichlet problem (\ref{dirichlet}).

It is easy to verify that $w+\tau\phi$ is a lower barrier for solutions of (\ref{dirichlet2}). Similarly, just using that
$H_K\ge 0$ we can see that  $-w+\tau\phi$ is a upper barrier for (\ref{dirichlet2}).
This concludes the proof for any value of $\tau$.
\qed

\section{Interior gradient estimates}

In this section, we establish an a priori global estimate for the
gradient.

\begin{proposition}\po\label{C1i}
Under the assumptions of Theorem \ref{main} there exists a positive
constant $C=C(\Omega,H,\phi, |\nabla z|_\Gamma|_0)$ independent of $\tau$ such that
$$
\sup_\Omega |\nabla z_\tau|\le C
$$
for any solution $z_\tau$ of the Dirichlet problem (\ref{dirichlet2}).
\end{proposition}

\noindent \emph{Proof:} The proof follows a similar
guideline as in \cite{DHL}. Consider on $\Sigma(z)$ the function
$$
\chi= ve^{2Kz},
$$
where $v=|\nabla z|^2 =z^i z_i$ and $K>0$ is a constant. We already have the desired bound if  $\chi$ achieves its maximum on $\Gamma$.
Thus, we assume that $\chi$ attains maximum value at an interior point $\bar u\in\Omega$ where $|\nabla z|> 0$. This assumption enables us to
choose a local normal coordinate system $x^{1},\ldots,x^{n}$ such that
$\partial_1|_{\bar u} = \nabla z/|\nabla z|$ and
$\sigma_{ij}(\bar u)=\delta_{ij}$. Hence, at
$\bar u$ we have
$$
z_1=|\nabla z|>0\;\;\mbox{and}\;\; z_j=0 \;\,\mbox{if}\;\, j\geq 2.
$$
Since $v_j = 2 z^l z_{l;j}$,
we also obtain at $\bar u$ from
$$
\chi_j = 2e^{2Kz}\big(Kvz_j + z^l z_{l;j}\big)=0
$$
that $z^l z_{l;j}=-Kvz_j$. Therefore, at  $\bar u$ we obtain
$$
z_{1;1}=-Kv, \quad v_1=-2Kv^{3/2} 
\quad\mbox{and}\quad v_j=0=z_{1;j}\;\,\mbox{if}\;\, j\ge 2.
$$
Moreover, we may assume  after rotating the $\partial_j$ that $(z_{i;j})$, $2\le i,j\le n$, is diagonal.

We write (\ref{basic2}) as 
\be\label{eqtilde} 
\tilde a^{ij}z_{i;j}=\tilde b, 
\ee 
where 
$$
\tilde a^{ij} =(\gamma+v)\sigma^{ij}-z^i z^j 
$$ 
and 
\be\label{tildeb} 
\tilde b=\frac{1}{2}\gamma_i z^i +
\big(\gamma+v\big)\left(\frac{\gamma_i}{2\gamma}z^i 
+ n\gamma\rho\right)+ n\lambda H (\gamma+v)^{3/2}. 
\ee
Hence,   at $\bar u$ we have  
\be\label{laplacian} 
(\gamma+v)z^i_{\,; i} = \tilde b - Kv^2. 
\ee  
Covariant differentiation of (\ref{eqtilde}) followed by
contraction with $\nabla z$  at $\bar u$ gives 
\be\label{eqder}
\frac{1}{\gamma+v}(\gamma_1v^{1/2}-2Kv^2)(\tilde b-Kv^2) -2K^2 v^3
+\tilde a^{ij}z^lz_{i;jl} = v^{1/2}\tilde b_{;1}. 
\ee
The third derivatives of $z$ satisfy the Ricci identity
$$
z_{i;jl}-z_{i;lj}=  R_{ikjl}\,z^k,
$$
where by $R_{ikjm}$ we denote the coefficients of the curvature tensor of
$M$.  Hence,
\begin{eqnarray*}
\tilde a^{ij}z^lz_{i;jl}
\!\!& =&\!\! \big( \gamma+v\big)\sigma^{ij}z^l\big(z_{l;ij}+
R_{ikjl}z^k\big) -z^i z^j z^l \big(z_{l;ij}+
R_{ikjl}z^k\big)\\
\!\!\!&=&\!\!\! \tilde a^{ij}z^lz_{l;ij}+(\gamma+v)R_{kl}z^kz^l\\
\!\!\!&=&\!\!\! \tilde a^{ij}z^lz_{l;ij}+R_{11}(\gamma+v)v
\end{eqnarray*}
since $R_{jklm}z^l z^m=0$. 

To estimate $\tilde a^{ij}z^lz_{l;ij}$ we use the Hessian matrix of $\chi$. Since
$$
\chi_{i;j}=2e^{2Kz}(2K^2z_i z_j v+ 2K z_i z^l z_{l;j}+Kz_{i;j} v
+2Kz_j z^{l}z_{l;i} + z^{l}_{\textrm{ };i}z_{l;j}+z^l z_{l;ij})
$$
is nonpositive at $\bar u$, it results that $\tilde
a^{ij}\chi_{i;j}\le 0$. Hence,
\begin{eqnarray*} 0 \!\!\!&\geq&\!\!\! 2K^2\tilde a^{ij}z_i z_j v+4K\tilde a^{ij}z_i
z^l z_{l;j}+K\tilde a^{ij}z_{i;j}v +\tilde a^{ij}z^{l}_{;i}z_{l;j}+\tilde a^{ij}z^lz_{l;ij}\\
\!\!\!&=&\!\!\!(\gamma+v)(Kvz^i_{\textrm{ };i} +z^l_{\textrm{
};i}z^i_{\textrm{ };l})-2K^2\gamma v^2
 +\tilde a^{ij}z^lz_{l;ij}.
\end{eqnarray*}
Since $(z_{i;k})$ is diagonal, it follows that
\be\label{ineq}
z^{l}_{\textrm{ };i}z^i_{\textrm{ };l}=(z^{l}_{\textrm{ };l})^2\geq
(z^{1}_{\textrm{ };1})^2= K^2v^2.
\ee
Using  (\ref{laplacian})  and  (\ref{ineq}) we obtain 
$$ 
\tilde a^{ij}z^l z_{l;ij}\le  K^2 \gamma v^2-K\tilde bv. 
$$ 
Thus,
\be\label{zlij}
\tilde a^{ij}z^lz_{i;jl}\le  K^2\gamma v^2-K\tilde bv+R(\gamma+v)v
\ee
where $R=R_{11}$.  Replacing (\ref{zlij}) into (\ref{eqder}) and multiplying
both sides  by $(\gamma+v)$ gives
\be\label{eq-der3} 
K\tilde b v
(\gamma+3v) + K^2\gamma v^2(v-\gamma)-Rv(\gamma+v)^2+ \gamma_1
v^{1/2}(Kv^2-\tilde b)+\tilde b_{;1} v^{1/2}(\gamma+v)\le 0. 
\ee
Since $\tilde b = \tilde b(x,z,\nabla z)$,  we have at $\bar u$ that
$$
\tilde b_{;1} =\tilde b_{x^1}+
\tilde b_{z}v^{1/2}-K \tilde b_{z^1}v.
$$
We rewrite (\ref{eq-der3}) as
\begin{eqnarray}
\label{eq-der4} & & K^2\gamma v^2(v-\gamma)-Rv(\gamma+v)^2+
\gamma_1v^{1/2}(Kv^2-\tilde b)
+\big(\gamma+v\big)\left(\tilde b_{z}v
+\tilde b_{x^1}v^{1/2}\right)\nonumber\\
& & \hspace*{14ex}+Kv\left(\tilde b (\gamma+3v)
-\tilde b_{z^1}(\gamma+v)v^{1/2}\right)\le 0.
\end{eqnarray}
From (\ref{tildeb}), we obtain
$$
\tilde b_{x^1} = 2n\rho\gamma\gamma_1
+\big(\gamma_{1;1}+n\gamma_1 \rho\big)v^{1/2}
+\frac{3}{2}n\lambda H \gamma_1 (\gamma+ v)^{1/2}+\big(\gamma_1/2\gamma\big)_{;1}v^{3/2}
+ n\lambda H_1 (\gamma+v)^{3/2}.
$$
Moreover,
\begin{eqnarray*}
\tilde b_{z} = n\gamma\rho_t(\gamma+v)+
n\lambda_t H (\gamma+v)^{3/2} \ge 0
\end{eqnarray*}
since $\rho_t\ge 0$ and $\lambda_t H\ge 0$.

We may assume from now on that $\gamma(\bar u)\le v(\bar
u)$. It follows that there exist constants $A_i=A_i(\gamma, \lambda, \rho, H, \nabla\gamma, \nabla H), \, 1\le i\le 4,$  such that
$$
\big(\gamma+v\big)\left(\tilde b_{z}v
+\tilde b_{x^1}v^{1/2}\right)
\ge\left(\left(\gamma_1/2\gamma\right)_{;1}
-2^{3/2}n\lambda|H_1|\right)v^3
+ A_1 v^2 + A_2 v^{3/2}+A_3 v + A_4 v^{1/2}.
$$
We also estimate
\begin{eqnarray*}
\gamma_1v^{1/2}(Kv^2-\tilde b)\ge K\gamma_1 v^{5/2}+B_1 v^2 
+ B_2v^{3/2}+B_3 v + B_4 v^{1/2},
\end{eqnarray*}
where $B_i=B_i(\gamma, \lambda, \rho, H, \nabla\gamma, \nabla H),\, 1\le i\le 4,$ are constants.
Finally, we have
$$
\tilde b (\gamma+3v)
-\tilde b_{z^1}(\gamma+v)v^{1/2}
=\gamma_1 v^{3/2}+ n\gamma \rho (\gamma+v)^2
+ n\lambda H \gamma (\gamma+v)^{3/2}\ge\gamma_1 v^{3/2}.
$$
With the above estimates, it follows from (\ref{eq-der4}) that
\begin{eqnarray*}
& &\left(K^2\gamma -R
+(\gamma_1/2\gamma)_{;1}-2^{3/2}n\lambda |H_1|\right)v^3 +
2K\gamma_1v^{5/2}\\
& &
\,\,\,\,\,\,\,-\big(K^2\gamma^2+2R\gamma-D_1\big)v^2+D_2v^{3/2}
-\big(R\gamma^2-D_3\big)v
+D_4v^{1/2} \le 0
\end{eqnarray*}
where $D_i=A_i+B_i$.
Choosing $K$ large enough so that the coefficient of $v^3$ is
positive implies that $v(\bar u)$ is bounded by a constant $K_0$
which depends only on $\Omega$ and $H$.

Using the fact that $\bar u$ is a maximum point for $\chi$, we
conclude that
$$
v(u) \leq K_0 e^{2K(z(\bar u)-z(u))}\le K_0 e^{4K|z|_0}.
$$
Finally, we observe that minor modifications in the calculations above give indeed
$$
|\nabla z_\tau | \le K_0^{1/2}e^{2K|z_0|}
$$
for any $\tau\in [0,1]$, and this concludes the proof. \qed

\section{Proof of corollaries}

In this section we prove Corollaries 2 and 3 in the Introduction.
\bigskip

\noindent \emph{Proof of Corollary \ref{main3}:}  We first observe  that the proofs of Propositions \ref{C1b} and \ref{C1i} still work under the weaker assumption $\inf_\Gamma
H_{K}\geq H\ge 0$. Thus, it suffices to show that Proposition \ref{C0} still works if $H^*\ge H$ when $Y$ is a Killing field. In this situation, it is easy to see that (\ref{dif}) is
equivalent to
$$
(H^*)^2(\gamma+ f'^2)^2f'^2-H^2(\gamma+ f'^2)^3
+2H^*\gamma D_0f'^2(\gamma+ f'^2) + \gamma^2 D_0^2f'^2>0.
$$
Clearly, the last  inequality holds for $D$ sufficiently large and the proof follows.\qed

\bigskip

\noindent \emph{Proof of Corollary \ref{main2}:} Being $Y$ closed we
may assume $\gamma=1$. Thus (\ref{kappa2}) and
(\ref{mean}) yield
$$
n^2H_K^2=(n-1)^2H_{\Gamma}^2.
$$
On the other hand,  the relation between the Ricci tensors of
$\bar M^{n+1}$ and $M^n$ is
$$
\mathrm{Ric}_{\bar M}(X,Z) ={\rm Ric}_M(X,Z)-(nk^2-k_t)\<X,Z\>
$$
for any $X,Z\in TM$. Thus (\ref{ric-amb}) is equivalent to
$$
\inf_{\Omega_0}{\mathrm{Ric}}_M^{rad}\ge -\frac{(n-1)^2}{n}\inf_{\Gamma}
H_\Gamma^2,
$$
and the proof follows.\qed\vspace{1ex}

\begin{remark}\po {\em We have from (\ref{kill-conf}) that the ambient space $\bar M^{n+1}$ is a product manifold  endowed with the metric
$$
\dd s^2=\lambda^2(t)\psi^2(u)(\dd t^2 + \psi^{-2}\dd\sigma^2).
$$
It is thus natural to consider the general situation of an ambient space $\R\times M^n$ endowed with  the  conformal metric
$$
\tilde g=\lambda^2(t,u) g=e^{2\varphi(t,u)} g
$$
where $g$ is the product metric in $\R\times M^n$. In this case,
the mean curvature equation for the graph $X=(z(u),u)$ is
$$
{\mathrm{div}}\, \bigg(\frac{\nabla z} {\sqrt{1+|\nabla
z|^2}}\bigg) -\frac{n}{\sqrt{1+|\nabla z|^2}}
\left(\langle \bar\nabla z,\bar\nabla \varphi\rangle -\varphi_t\right)-n\lambda H=0,
$$
where $z(t,u)=z(u)$ and we compute $\langle \bar\nabla z,\bar\nabla \varphi\rangle$ in the ambient space.
We easily conclude  that in order to have a maximum principle for the above equation we have to ask $\varphi_t$ to be independent of $t$, that is, the function  $\lambda$ has to separate variables. But this is  precisely the case we studied in this paper.
}\end{remark}

{\renewcommand{\baselinestretch}{1} \hspace*{-20ex}\begin{tabbing}
\indent \= Marcos Dajczer\\
\> IMPA \\
\> Estrada Dona Castorina, 110\\
\> 22460-320 -- Rio de Janeiro -- Brazil\\
\> marcos@impa.br\\
\end{tabbing}}

\vspace*{-6ex}

{\renewcommand{\baselinestretch}{1} \hspace*{-20ex}\begin{tabbing}
\indent \= Jorge Herbert S. de Lira\\
\> UFC - Departamento de Matematica \\
\> Bloco 914 -- Campus do Pici\\
\> 60455-760 -- Fortaleza -- Ceara -- Brazil\\
\> jorge.lira@pq.cnpq.br
\end{tabbing}}

\end{document}